\documentstyle[11pt]{article}
\input{amssym.def}
\input epsf

\newtheorem{thm}{Theorem}[section]
\newtheorem{lemma}[thm]{Lemma}
\newtheorem{prop}[thm]{Proposition}
\newtheorem{conj}[thm]{Conjecture}

\newtheorem{cor}[thm]{Corollary}
\newtheorem{rmk}[thm]{Remark}
\newtheorem{ex}[thm]{Example}

\newenvironment{pf}{\bf Proof:\rm}{\begin{flushright}$\Box$\end{flushright}}
\newcommand{\Qed}{\begin{flushright}$\Box$\end{flushright}}
\newcommand{\as}{{\rm{as}}}

\def\rank{\mathop{\rm rank}\nolimits}

\def\mapdown#1{\Big\downarrow\rlap{$\vcenter{\hbox{$\scriptstyle#1$}}$}}

\newcommand{\C}{{\Bbb C}}
\newcommand{\D}{{\rm D}}
\newcommand{\Z}{{\Bbb Z}}
\newcommand{\R}{{\Bbb R}}

\newcommand{\HH}{{\Bbb H}}
\newcommand{\cS}{{\cal S}}

\newcommand{\cB}{{\cal B}}

\newcommand{\cH}{{\cal H}}
\newcommand{\cL}{{\cal L}}
\newcommand{\cN}{{\cal N}}
\newcommand{\cO}{{\cal O}}
\newcommand{\cP}{{\cal P}}

\newcommand{\cU}{{\cal U}}

\newcommand{\cX}{{\cal X}}

\newcommand{\pF}{{\cal F}}
\newcommand{\IC}{{\cal IC}}
\newcommand{\pj}{{{{}^{\scriptscriptstyle p}} \hspace{-.02in} j}}
\newcommand{\fg}{{\frak g}}
\newcommand{\ft}{{\frak t}}
\newcommand{\fb}{{\frak b}}
\newcommand{\bb}{{\backslash\!\backslash}}

\begin{document}
\bibliographystyle{plain}

\title{A Generalization of Springer Theory Using Nearby Cycles}
\author{Mikhail Grinberg}
\date{May 20, 1998}
\maketitle

\begin{abstract}
Let $\fg$ be a complex semisimple Lie algebra, and $f : \fg \to G \bb \fg$
the adjoint quotient map.  Springer theory of Weyl group representations
can be seen as the study of the singularities of $f$.

In this paper, we give a generalization of Springer theory to visible,
polar representations.  It is a class of rational representations of
reductive groups over $\C$, for which the invariant theory works by
analogy with the adjoint representations.  Let $G \, | \, V$ be such a
representation, $f : V \to G \bb V$ the quotient map, and $P$ the sheaf
of nearby cycles of $f$.  We show that the Fourier transform of $P$ is
an intersection homology sheaf on $V^*$.

Associated to $G \, | \, V$, there is a finite complex reflection group
$W$, called the Weyl group of $G \, | \, V$.  We describe the endomorphism
ring $\rm{End} (P)$ as a deformation of the group algebra $\C [W]$.
\end{abstract}

\section{Introduction}

Let $f: \C^d \to \C$ be a non-constant polynomial.  Fix a point
$x \in E = f^{-1} (0)$.  The {\em Milnor fiber} of $f$ at $x$ is
defined by $F_{f,x} = f^{-1} (\epsilon) \cap {\rm B}_{x,\delta}$,
where ${\rm B}_{x,\delta}$ is the $\delta$-ball around $x$, and
we assume that $0 < |\epsilon| \ll \delta \ll 1$.  Since the
publication of Milnor's seminal work [Mi], the cohomology of
$F_{f,x}$ has been a central object in the study of the
singularities of $f$.

In the language of sheaf theory, the cohomology groups
$H^* (F_{f,x})$, for all $x \in E$, fit together to form a
constructible, bounded complex $P$ of sheaves on $E$, called
the {\em sheaf of nearby cycles} of $f$ (i.e., $H^* (F_{f,x})$
is the stalk of $P$ at $x$).  The complex $P$ is a perverse
sheaf ([GM2], [KS]).

There are two kinds of questions one may ask about the
singularities of $f$.  First, one may be interested in the
cohomology of the Milnor fibers $F_{f,x}$ and of other related
spaces, as well as in the various maps between such groups.
The most important among such maps is the monodromy
transformation $\mu_x : H^* (F_{f,x}) \to H^* (F_{f,x})$,
arising from the dependence of the Milnor fiber on the choice
of the small number $\epsilon \in \C^*$.  Much of the work in
singularity theory since the publication of [Mi] has centered
around such concrete geometric questions.

On the other hand, one may ask: what is the structure of $P$
as an object in the abelian, artinian category of perverse
sheaves on $E$ (see [BBD])?  In particular, what are the simple
constituents of $P$?  Is $P$ semisimple?  What are the extensions
involved in building $P$ up from its simple constituents?  What
is the endomorphism ring $\rm{End} (P)$?  How does the monodromy
transformation $\mu \in \rm{Aut} (P)$ act?

These two kinds of questions are intimately related.  The
structure of $P$ as a perverse sheaf is the ``glue'' that ties
together the various local geometric invariants of $f$ at different
points.  To the author's knowledge, the only general theorem about
the structure of $P$ is the deep result of Ofer Gabber (see [BB])
describing the filtration of $P$ arising from the action of the
monodromy.  Gabber's theorem is a semisimplicity assertion about
certain subquotients of $P$.

The definition of the nearby cycles sheaf $P$ can be extended to
the setting of a dominant map $f : \C^d \to \C^r$.  In order to
assure the desired properties of $P$, one needs to impose some
technical conditions on $f$ (see Section 2.2 below).  On the
other hand, the resulting theory is richer, because instead of
a single monodromy transformation of $P$, we have an action
of the whole fundamental group of the set of regular values of
$f$ near zero.

In this paper, we consider the nearby cycles for a certain special
class of maps $f : \C^d \to \C^r$, whose components are given by
homogeneous polynomials.  This class arises naturally from the
invariant theory of reductive group actions, and our motivation
for studying it comes from the Springer theory of Weyl group
representations.

Springer theory (see [BM1], [M], [S3], and Section 2.1 below) can
be seen as the study of the singularities of the adjoint quotient
map $f: \fg \to G \bb \fg = \rm{Spec} \, \C[\fg^*]_{}^G$,
associated to a complex semisimple Lie algebra.  For this map,
both the categorical structure of $P$ and its local invariants are
well understood.  The sheaf $P$ is semisimple.  The monodromy group
acting on $P$ is the Braid group $B^{}_W$ of the Weyl group $W$ of
$\fg$; but the monodromy action factors through $W$.  Moreover,
this monodromy action gives an isomorphism $\C[W] \cong \rm{End}
(P)$.  Thus, Springer theory gives a relation between the
singularities of the nilcone $\cN = f^{-1} (0) \subset \fg$ and the
representations of the Weyl group.

Dadok and Kac [DK] introduced a class of rational representations
$G \, | \,V$ ($G$ reductive / $\C$) for which the invariant theory
works by analogy with the adjoint representations.  They call this
the class of {\em polar representations}; it includes the adjoint
representations, as well as many of the classical invariant
problems of linear algebra (see Section 3.2 for examples).  For any
polar representation $G \, | \, V$, the quotient $G \bb V$ is
isomorphic to a vector space.  In this paper, we study the
singularities of the quotient map $f : V \to G \bb V$, giving a
generalization of Springer theory to polar representations
satisfying a mild additional hypothesis.

Our main result (Theorem \ref{main}, part (i)) is the following.
For a polar representation $G \, | \, V$, assume that the fiber
$E = f^{-1} (0)$ consists of finitely many $G$-orbits (this
condition is called {\em visible}).  Then the nearby cycles sheaf
$P$ of the quotient map $f$ satisfies:
\begin{equation}\label{intro}
\pF \, P \cong \IC ((V^*)^{rs}, \cL),
\end{equation}
where $\pF$ is the geometric Fourier transform functor, and the
right hand side is an intersection homology sheaf on $V^*$.
(See [KS] for a definition of the Fourier transform, and [GM1]
for a discussion of intersection homology.)  The content of this
is that the sheaf $P$ is completely encoded in the single local
system $\cL$ on a certain locus in the dual space $V^*$.  In the
case of Springer theory, this result is due independently to
Ginzburg [Gi] and to Hotta-Kashiwara [HK] (see also [Br]).
The proof of (\ref{intro}) draws on a general result of the
author about the specialization of an affine variety to the
asymptotic cone (see [Gr1] and Section 2.2.2).

Further assertions of Theorem \ref{main} describe the holonomy
of the local system $\cL$ and the monodromy action on $P$ of
the appropriate fundamental group.  It turns out that the
semisimplicity observed in Springer theory does not extend to
this generalization.  Associated to each polar representation
$G \, | \, V$, there is a finite complex reflection group $W$,
called the Weyl group of $G \, | \, V$ (it is not, in general,
the Weyl group of $G$).  We have $\dim \, \rm{End} (P) =
|W|$, but the algebra $\rm{End} (P)$ is not, in general,
isomorphic to $\C [W]$.  Instead, it is given as a kind of a
Hecke algebra associated to $W$.

We should note that our Fourier transform description of $P$,
while well suited to the study of the endomorphism ring and of
the monodromy action, falls short of giving the complete structure
of $P$ as a perverse sheaf.  This is because intersection homology
is not an exact functor from local systems to perverse sheaves.

An important special case of this theory, which includes the
adjoint representations, arises in the following way.  Let
$\theta : \fg \to \fg$ be an involutive automorphism of a complex
semisimple Lie algebra, and $\fg = \fg_+ \oplus \fg_-$ the
eigenspace decomposition for $\theta$.  Then the adjoint form
$G_+$ of the Lie algebra $\fg_+$ acts on the symmetric space
$\fg_-$.  Orbits and invariants of the representation
$G_+ \, | \, \fg_-$ were studied by Kostant and Rallis in [KR].
This representation is polar and visible.  In Theorem \ref{symsp},
we explicitly compute the endomorphism ring $\rm{End} (P)$ in this
case.  It is given as a ``hybrid'' of the group algebra $\C[W]$
and the Hecke algebra of $W$ (which is a Coxeter group)
specialized at $q = - 1$.

I would like to thank the University of Utrecht, Holland, and the
IAS, Princeton, for their hospitality.  Throughout this work I was
supported by a Fannie and John Hertz foundation graduate fellowship.
Discussions with Tom Braden, Ian Grojnowski, Victor Kac, Peter
Magyar, David Massey, Dirk Siersma, and Tony Springer have been
of great value to me.  I am also grateful to Mark Goresky, David
Kazhdan, and Wilfried Schmid for their interest and encouragement.
Finally, I would like to thank my adviser, Robert MacPherson, for
his inspiration, guidance, and support.

\tableofcontents

\pagebreak

\section{Background}

{\bf Notations.}
We will say {\em sheaf} to mean {\em complex of sheaves}
throughout; all our sheaves will be sheaves of $\C$ vector
spaces.  Given a map $g : X \to Y$, the symbols $g_*$, $g_!$
will always denote the {\em derived} push-forward functors.
All perverse sheaves and intersection homology will be taken
with respect to the middle perversity (see [GM1], [BBD]);
we use the shift conventions of [BBD].  Given a sheaf $A$ on
$X$, and a pair of closed subspaces $Z \subset Y \subset X$,
we will write $\HH^k(Y, Z; A)$ for the hypercohomology group
$\HH^k(j_! \,i^* A)$, where $i : Y \setminus Z \to X$ and
$j : Y \setminus Z \to Y$ are the inclusion maps.  We call
$\HH^k(Y, Z; A)$ the relative hypercohomology of $A$.

For an analytic function $f : M \to \C$, we use the notation of
[KS, Chapter 8.6] for the nearby and vanishing cycles
functors $\psi^{}_f$, $\phi^{}_f$.  Let $A$ be a sheaf on $M$
or $M \setminus f^{-1} (0)$.  We denote by $\mu : \psi^{}_f A \to
\psi^{}_f A$ the monodromy transformation of the nearby cycles. 
It is the counter-clockwise monodromy for the family of sheaves
$\psi^{}_{f / \tau} A$, parametrized by the circle
$\{ \tau \in \C \, | \, | \tau | = 1 \}$.

When $V$ is a $\C$ vector space, we denote by $\cP^{}_{\C^*} \,
(V)$ the category of $\C^*$-conic perverse sheaves on $V$, and
by $\pF : \cP^{}_{\C^*} \, (V) \to \cP^{}_{\C^*} \, (V^*)$
the (shifted) Fourier transform functor.  In the notation
of [KS, Chapter 3.7], we have $\pF \, P = P \, \hat{} \;
[\dim \, V]$.  To avoid cumbersome notation, we will use
the following shorthand: if $P$ is a conic perverse sheaf on a
closed, conic subvariety $X \subset V$, and $j^{}_X : X \to V$
is the inclusion, then we write $\pF \, P$ instead of
$\pF \circ (j^{}_X)_* \, P$.  When $E$ is a $G$-space, we denote
by $\cP^{}_G \, (E)$ the category of $G$-equivariant perverse
sheaves on $E$.  The symbol $\dim$, without a subscript, will
always denote the {\em complex} dimension.

\subsection{Springer Theory}

In this section we give a brief summary of the Springer theory
of Weyl group representations.  Springer theory (see [BM1], [S3])
is concerned with exhibiting the Weyl group $W$ of a complex
semisimple Lie algebra $\fg$ as the symmetry group of a certain
perverse sheaf $P$ on the nilcone $\cN \subset \fg$.  As we
mentioned in the introduction, $P$ is the nearby cycles sheaf
for the adjoint quotient map $f : \fg \to G \bb \fg$ (here
$G$ is the adjoint form of $\fg$).  However, this is not the
original definition of $P$, and not the one used in the
literature to prove its properties.

We now describe the construction due to Lusztig [Lu1] and
Borho-MacPherson [BM1] of the sheaf $P$ and the $W$-action on
it.  Consider the Grothendieck simultaneous resolution diagram:

\vspace{.1in}

$$\def\normalbaselines{\baselineskip20pt
\lineskip3pt \lineskiplimit3pt}
\def\incmapright#1{\smash{
\mathop{\hookrightarrow}\limits^{#1}}}
\def\mapright#1{\smash{
\mathop{\longrightarrow}\limits^{#1}}}
\def\mapdown#1{\Big\downarrow\rlap
{$\vcenter{\hbox{$\scriptstyle#1$}}$}}
\matrix{\tilde{\cN}  & \incmapright{}   & \tilde{\fg} &
\mapright{\tilde{f}} & \ft              \cr
        \mapdown{p}  &                  & \mapdown{q} &
                     & \mapdown{f|_\ft} \cr
        {\cN}        & \incmapright{j}  & \fg         &
\mapright{f}         & W \, \bb \, \ft  \cr}$$

\vspace{.2in}

\noindent
Here $\tilde{\fg}$ is the variety of pairs $(x, \fb)$, where
$\fb \subset \fg$ is a Borel subalgebra, and $x \in \fb$.  This
variety is smooth; it is a vector bundle over the flag variety
$\cB$ of $\fg$.  The map $f$ is given by the invariants of the
adjoint action $G \, | \, \fg$; the target space is identified
with the quotient of a Cartan subalgebra $\ft \subset \fg$ by
the action of the Weyl group.  The map $\tilde{f}$ is the
natural morphism making the diagram commute.  Finally,
$\tilde{\cN} = \{ \, (x, \fb) \in \tilde{\fg} \; | \;
x \in \cN \, \}$.

Note that the map $q$ is a finite cover over the set $\fg^{rs}$
of regular semisimple elements in $\fg$.  Given a point
$x \in \fg^{rs}$, the fiber $q^{-1}(x)$ can be identified with
the Weyl group $W$.  Such an identification requires a choice of
a positive Weyl chamber for the Cartan subalgebra
$\fg_x \subset\fg$.  The Weyl group $W$, therefore, acts as the
deck transformations of the covering
$q^{rs} : \tilde{\fg}^{rs} \to \fg^{rs}$ (where
$\tilde{\fg}^{rs} = q^{-1} (\fg^{rs})$, and $q^{rs}$ is the
restriction of $q$).

The diagram above is a simultaneous resolution in the following
sense.  For any point $t \in \ft$, the restriction
$q^{}_t : \tilde{f}^{-1} (t) \to f^{-1} (f(t))$ is a resolution
of singularities.  In particular, the map $p = q^{}_0$ is a
resolution of singularities of the nilcone.

We now consider the push-forward sheaves
$$P = p^{}_* \, \C_{\tilde{\cN}} \, [\dim \, \cN] \;\;
\mbox{and} \;\;
Q = q^{}_* \, \C_{\tilde{\fg}} \, [\dim \, \fg].$$
Note that $P = j^* \, Q \, [-r]$, where $r =  \dim \, \ft$.

\begin{prop}\label{dimfib} {\em [Lu1]}
The map $q$ is small and the map $p$ is semi-small in the sense
of {\em [GM1]}.  Consequently, the sheaf $P$ is perverse, and
the sheaf $Q$ is the intersection cohomology extension
$\IC(\fg^{rs}, \cL)$ of the local system
$\cL = q^{rs}_* \, \C_{\tilde{\fg}^{rs}}$ on $\fg^{rs}$.
\end{prop}

The Weyl group action on $\tilde\fg^{rs}$ by deck
transformations produces an action on the local system $\cL$,
and by the functoriality of intersection cohomology on the
sheaf $Q$.  In fact, we have ${\rm End} (\cL) = {\rm End} (Q)
= \C [W]$, the group algebra of $W$.  Using the relation
$P = j^* \, Q \, [-r]$, we obtain a $W$-action on $P$ as well.
This construction of the $W$-action is due to Lusztig [Lu1].
Originally, the representations of $W$ on the stalks of $P$,
i.e., on the cohomology of the varieties $p^{-1} (x)$
($x \in \cN$), were constructed by Springer ([S1], [S2]).
Other constructions of the same representations were given by
Slodowy [Sl] and Kazhdan-Lusztig [KL].  The construction of
Slodowy is essentially the nearby cycles definition, which is
the starting point for this paper.  Hotta in [Ho] showed that
all of the different constructions agree.

\begin{thm}\label{bmmain} {\em [BM1]}

(i)  $P$ is a semisimple perverse sheaf.

(ii) The action of $W$ on $P$ gives an isomorphism
$\C [W] \cong {\rm End} (P)$.
\end{thm}

Borho and MacPherson deduce the first part of Theorem
\ref{bmmain} from the decomposition theorem of Beilinson,
Bernstein, Deligne, and Gabber [BBD], and use a counting
argument to prove the second part.  A different proof, using
the Fourier transform, was later given by Ginzburg [Gi] and
Hotta-Kashiwara [HK] (see also [Br]).  We will use the Killing
form on $\fg$ to identify $\fg$ with its dual, and to regard the
Fourier transform on $\fg$ as a functor $\pF : \cP^{}_{\C^*}
(\fg) \to \cP^{}_{\C^*} (\fg)$.

\begin{thm}\label{ftinspth} {\em [Gi], [HK]}

(i)  There is an isomorphism $\pF \, P \cong Q$.

(ii) The Weyl group $W$ acts on each side of this isomorphism.
These two actions differ by the sign character of $W$ (i.e.,
the character that sends each simple reflection $\sigma \in W$
to $-1$.)
\end{thm}

Theorem \ref{ftinspth} implies Theorem \ref{bmmain}, because
the functor $\pF$ is an equivalence of categories.

\subsection{Nearby Cycles}

\subsubsection{Families Over a Large Base}

The sheaf $P$ constructed in the previous section can be defined in
a different way, namely, as the nearby cycles of the adjoint quotient
map $f$.  One usually considers the nearby cycles functor for an
analytic function $M \to \C$ (see [GM2, Section 6], [KS, Chapter 8.6]).
Proposition \ref{nclb} below provides a technical basis for discussing
the nearby cycles for a sufficiently nice map $M \to \C^r$.

Let $M$ be a connected complex analytic manifold of dimension $d$,
$\cU$ a neighborhood of zero in $\C^r$, and $f: M \to \cU$ an analytic
map onto $\cU$.  Let $\cU^{sing} \subset \cU$ be the closure of the
set of non-regular values of $f$.  We assume that $\cU^{sing}$ is a
proper analytic subvariety of $\cU$.  Let $\cU^{reg} =
\cU \setminus \cU^{sing}$, and note that the preimage
$M^\circ = f^{-1} (\cU^{reg}) \subset M$ is a manifold.  We assume
that there is an analytic Whitney stratification $\cS$ of the fiber
$E = f^{-1} (0)$, such that Thom's ${\rm A}_f$ condition holds for
the pair $(M^\circ, S)$, for every stratum $S \in \cS$.  Recall that
the ${\rm A}_f$ condition says that for any sequence of points
$x_i \subset M^\circ$, converging to a point $y \in S$, if there
exists a limit
$$\Delta = \lim_{i \to \infty} T_{x_i} f^{-1} (f (x_i))
\subset T_y M,$$
then $\Delta \supset T_y S$.  If $r = 1$, such a stratification
$\cS$ can always be found.  For $r > 1$, the existence of $\cS$
is an actual restriction on $f$.  It implies, in particular,
that $\dim E = d - r$.  We refer the reader to [Hi] for a detailed
discussion of the ${\rm A}_f$ condition.

Let $U$ be a neighborhood of zero in $\C$, and $\gamma : U \to \cU$
an embedded analytic arc, such that $\gamma (0) = 0$, and
$\gamma (z) \in \cU^{reg}$ for $z \neq 0$.  Form the fiber product
$M_\gamma =  M \times^{}_{\cU} U$, and let $f_\gamma : M_\gamma
\to U$ be the projection map.  We may consider the nearby cycles
sheaf $P_{f_\gamma} = \psi_{f_\gamma} \, \C_{M_\gamma} \, [d - r]$.

\begin{prop}\label{nclb}
(i)   The sheaves $P_{f_\gamma}$ for different $\gamma$ are all
isomorphic.  We may therefore omit the subscript $\gamma$, and
call the sheaf $P_f = P_{f_\gamma}$ the nearby cycles of $f$.
It is a perverse sheaf on $E$, constructible with respect to $\cS$.

(ii)  The local fundamental group $\pi_1 (\cU^{reg} \cap
{\rm B}_\epsilon)$, where ${\rm B}_\epsilon \subset \cU$
is a small ball around the origin, acts on $P_f$ by monodromy.  We
denote this action by $\mu : \pi_1 (\cU^{reg} \cap
{\rm B}_\epsilon) \to {\rm Aut} (P_f)$.

(iii) The sheaf $P_f$ is Verdier self-dual.
\end{prop}
\begin{pf}
We refer the reader to [GM2, Section 6] for a proof that each
$P_{f_\gamma}$ is a perverse sheaf constructible with respect
to $\cS$.  The self-duality of $P_{f_\gamma}$ follows from the 
fact that the nearby cycles functor $\psi_{f_\gamma}$ commutes
with Verdier duality.

For the other assertions of the proposition, choose a point
$x \in E$, and let $S \subset E$ be the stratum containing $x$.  Fix
a normal slice $N \subset M$ to $S$ through the point $x$.  Also fix
a pair of numbers $0 < \epsilon \ll \delta \ll 1$ (chosen to be small
in decreasing order).  For each regular value $\lambda \in \cU^{reg}$
with $|\!| \lambda |\!| < \epsilon$, consider the Milnor fiber
cohomology group
$$H^* (f^{-1} (\lambda) \cap N \cap {\rm B}_{x, \delta}),$$
where ${\rm B}_{x, \delta} \subset M$ is the $\delta$-ball around $x$
(we fix a Hermitian metric on $M$ for this).  It is enough to show
that this Milnor fiber cohomology varies as a local system in
$\lambda$.  This can be done by adapting the argument in
[L\^e1, Section 1].
\end{pf}

Returning to the notation of the previous section, it is easy to
check that the adjoint quotient map $f : \fg \to W \, \bb \, \ft$
and the orbit stratification of the nilcone $\cN$ satisfy the
hypotheses of Proposition \ref{nclb}.  Note that the set of regular
values of $f$ equals $f (\ft^{reg})$, where $\ft^{reg} =
\ft \cap \fg^{rs}$.  The fundamental group $\pi_1 (f (\ft^{reg})
\cap {\rm B}_\epsilon)$ is the braid group $B^{}_W$ associated to $W$.
Theorem \ref{ncinspth} below is a consequence of the properties of the
simultaneous resolution (see [Sl] for an early treatment of these
ideas).

\begin{thm}\label{ncinspth} {\em [Ho], [M]}

(i)  We have $P \cong P_f$.

(ii) The monodromy action $\mu : B^{}_W \to {\rm Aut} (P_f)$ factors
through the natural homomorphism $B^{}_W \to W$, producing an action
of $W$ on $P_f$.  The isomorphism of part (i) agrees with the $W$
actions on both sides.
\end{thm}

Thus, Springer theory can be viewed as the study of the nearby cycles
of the adjoint quotient map.

\pagebreak

\subsubsection{Specialization to the Asymptotic Cone}

In this section we summarize the results of [Gr1] on which this
paper relies.

Let $V \cong \C^d$ be a complex vector space, and $X \subset V$
a connected, smooth, closed subvariety.  We denote by $\bar V$
the standard projective compactification of $V$, and by $\bar X$
the closure of $X$ in $\bar V$.  Set $V^\infty = \bar V
\setminus V$, and $X^\infty = \bar X \cap V^\infty$.  The
asymptotic cone $\as (X) \subset V$ is defined as the affine
cone over $X^\infty$.

Another way to define $\as (X)$ is as follows.  Let 
$\tilde X^\circ = \{ (\lambda, \tilde x) \in \C^* \times V \, |
\, \tilde x \in \lambda X \}$, and $\tilde X$ be the closure
of $\tilde X^\circ$ in $\C \times V$.  Write $\tilde f :
\tilde X \to \C$ for the projection on the first factor.
Then $\as (X) = \tilde f^{-1} (0)$.  This definition shows that
the asymptotic cone is naturally equipped with a nearby cycles
sheaf $P = P^{}_X := \psi_{\tilde f} \, \C_{\tilde X} [n]$, where
$n = \dim X$.  The sheaf $P$ is $\C^*$-conic.

Given a Hermitian inner product on $V$, we may consider for any
$x \in F$, the angle $\angle (x, T_x X) \in [0, \pi/2]$ between
the vector $x \in V$ and the subspace $T_x X \subset T_x V \cong V$.
The variety $X$ is said to be {\em transverse to infinity} if for
some (equivalently for any) inner product on $V$, there exists
a constant $k > 0$, such that for any $x \in X$, we have:
$$\angle (x, T_x X) < \frac{k}{|\!| x |\!|} \, .$$

\begin{thm}\label{asym_cone}{\em [Gr1, Theorem 1.1]}
Assume $X \subset V$ is transverse to infinity and let $P$ be
the nearby cycles sheaf on $\as (X)$.  Let $T_X^* V \subset T^* V
\cong V \times V^*$ be the conormal bundle to $X$, and
$p_2 : V \times V^* \to V^*$ be the projection.  Let
$Y = \overline{p_2 (T_X^* V)}$; it is an irreducible cone in $V^*$.
Then we have:
$$\pF \, P \cong \IC (Y^\circ, \cL),$$
where $\pF$ is the geometric Fourier transform functor, $\cL$ is
a local system on some Zariski open subset of $Y^\circ \subset Y$,
and the right hand side is the intersection homology extension of
$\cL$. 
\end{thm}

We will need some auxiliary facts describing the stalks of the
$\pF \, P$.

\begin{lemma}\label{odlim}{\em [Gr1, Proposition 3.3]}
In the situation of Theorem \ref{asym_cone}, fix $l \in V^*$,
and let $\xi = \mbox{\em Re} \, (l)$.  Also fix large positive
numbers $1 \ll \xi^{}_0 \ll \eta^{}_0$.  Then we have:
$$H^i_l \, (\pF \, P) \cong H^{i+d+n} \,
( \{ \, x \in X \, | \, \xi (x) \leq \xi^{}_0 \, \}, \,
\{ \, |\!| x |\!| \geq \eta^{}_0 \, \}; \, \C ).$$
\end{lemma}

Assume now $l \neq 0$.  Let $\Delta \subset V$ be the kernel
of $l$, and $L \subset V$ be any line complementary to $\Delta$.
We have $V = \Delta \oplus L$.  Take the standard projective
compactification $\bar \Delta$ of $\Delta$, and let $\hat V  =
\hat V_l = \bar \Delta \times L$.  It is not hard to check that
the space $\hat V$ is canonically independent of the choice of the
line $L$.  Note that $l : V \to \C$ extends to a proper algebraic
function $\hat l : \hat V \to \C$.  Let $\hat X = \hat X_l$ be the
closure of $X$ in $\hat V$, and $j : X \to \hat X$ be the inclusion
map.  Set $\hat X^\infty = \hat X \setminus X$.

\begin{lemma}\label{auxmtw}{\em [Gr1, Proposition 3.7]}
For $l \neq 0$, the statement of Lemma \ref{odlim} can be modified
as follows:
$$H^i_l \, (\pF \, P) \cong \HH^{i+d} \, ( \{ \, x \in \hat X
\, | \, \hat \xi (x) \leq \xi^{}_0 \, \}, \, \{ \, |\!| \hat l
(x) |\!| \geq 2 \xi^{}_0 \, \}; \; j^{}_! \, \C^{}_X \, [n] ),$$
where $\xi^{}_0 \gg 1$, and $\hat \xi (x) = \mbox{\em Re} \,
(\hat l(x))$.
\end{lemma}

Let $\cS$ be an algebraic Whitney stratification of $\as (X)$,
written $\as (X) = \bigcup_{S \in \cS} S$, satisfying the following
three conditions.

(i)   $\cS$ is conic, i.e., each $S \in \cS$ is $\C^*$-invariant.

(ii)  Thom's ${\rm A}_{\tilde f}$ condition holds for the pair
$(\tilde X^\circ, S)$, for each $S \in \cS$.

(iii) Let $\cS^\circ = \cS \setminus \{\{ 0 \}\}$.  For
$S \in \cS^\circ$, let $S^\infty \subset X^\infty$ be the
projectivization of $S$.  Then the decomposition $\bar X =
X \cup \bigcup_{S \in \cS^\circ} S^\infty$ is a Whitney
stratification.

The existence of such an $\cS$ follows form the general results
of stratification theory.  Let $a (\cS, l)$ be the dimension of
the critical locus of $l \, |_{\as (X)}$ with respect to $\cS$.

\begin{lemma}\label{dimbound}{\em [Gr1, Section 3.4]}
There exists a stratification $\hat \cX$ of $\hat X$ with the
following property.  Write $Z$ for the critical locus of
$\hat l \, |_{\hat X}$ with respect to $\hat \cX$.  Then
$$\dim \, Z \cap \hat X^\infty < a (\cS, l).$$
(We set $\dim \, \emptyset = -1$.)
\end{lemma}

\subsection{Polar Representations}

\subsubsection{A Summary of the Results of Dadok and Kac}

Dadok and Kac [DK] introduced and studied the class of polar
representations.  Their motivation was to find a class of
representations of reductive groups over $\C$, whose
invariant theory works by analogy with the adjoint
representations.  In this section we will, following [DK],
recall the definition and the main properties of polar
representations.

Let $G \, | \, V$ be a rational linear representation of a
connected reductive Lie group over $\C$ in a finite dimensional
vector space.  A vector $v \in V$ is called {\em semisimple} if
the orbit $G \cdot v$ is closed, and {\em nilpotent} if the
closure $\overline{G \cdot v}$ contains zero.  We say that
$v \in V$ is {\em regular semisimple} if $G \cdot v$ is closed
and of maximal dimension among all closed orbits.  We will write
$V^s$ ($V^{rs}$) for the set of all semisimple (regular semisimple)
vectors in $V$.  The representation $G \, | \, V$ is called
{\em stable} if $\dim \, G \cdot v \geq \dim \, G \cdot x$, for
any regular semisimple $v \in V$ and any $x \in V$.

For a semisimple vector $v \in V$, define a subspace
$$c_v = \{ \, x \in V \; | \; \fg \cdot x \subset
\fg \cdot v \, \},$$
where $\fg$ is the Lie algebra of $G$.  The orbits through
$c_v$ thus have ``parallel'' tangent spaces.  The representation
$G \, | \, V$ is called {\em polar} if for some semisimple
$v \in V$, we have $\dim \, c_v = \dim \, \C [V^*]_{}^G$;
in this case, $c_v$ is called a {\em Cartan subspace}.
The prototype of all polar representations (responsible for the
name ``polar'') is the action of the circle $S^1$ on the plane
$\R^2$ by rotations.  The class of polar representations
includes the adjoint representations (a Cartan subspace for an
adjoint representation is just a Cartan subalgebra) and the
representations arising from symmetric spaces, studied by Kostant
and Rallis in [KR] (see Section 6 below).  It is, in fact, much
larger.  Every representation with $\dim \, \C [V^*]_{}^G = 1$ is
automatically polar.  Dadok and Kac [DK] give a complete list of
all polar representations of simple $G$.  Their classification
includes many of the classical invariant problems of linear algebra.
Other examples of polar representations can be found in [K].  Note
that a polar representation need not be stable.

Before describing the invariant theory of a polar representation,
we recall the notion of a Shephard-Todd group (or a complex
reflection group).  Let $c$ be a complex vector space.  An element
$g \in GL(c)$ is called a {\em complex reflection} if in some
basis it is given by the matrix $\mbox{diag} \, (e^{2 \pi i / n},
1, \dots , 1)$, for some integer $n > 1$.  A finite subgroup
$W \subset GL(c)$ is called a {\em Shephard-Todd group} if it is
generated by complex reflections.  Shephard-Todd groups are a
natural generalization of Coxeter groups; we refer the reader to
[ST] for a discussion of their properties.  Let $W \subset GL(c)$
be a Shephard-Todd group, $g$ a complex reflection in $W$.  Fix a
basis of $c$ in which $g = \mbox{diag} \, (e^{2 \pi i / n}, 1,
\dots , 1)$.  We will say that $g$ is a {\em primitive reflection}
if there is no integer $n' > n$ such that $\mbox{diag} \,
(e^{2 \pi i / n'}, 1, \dots , 1) \in W$.

\begin{thm}\label{polar} {\em [DK]}
Let $G \, | \, V$ be a polar representation of a connected
reductive group.

(i)    All Cartan subspaces of $V$ are $G$-conjugate.

(ii)   Fix a Cartan subspace $c \subset V$.  Every vector
$v \in c$ is semisimple and every closed orbit passes
through $c$.

(iii)  Let $N^{}_G (c)$ and $Z^{}_G (c)$ be, respectively, the
normalizer and the stabilizer of $c$.  The quotient $W = N^{}_G (c)/
Z^{}_G (c)$ is a Shephard-Todd group acting on $c$; it is called the
{\em Weyl group} of the representation $G \, | \, V$ (note that, in
general, $W$ is not the Weyl group of $G$).

(iv)    Restriction to $c$ gives an isomorphism of invariant
rings $\C [V^*]_{}^G \cong \C [c^*]_{}^W$; this invariant ring
is free, generated by homogeneous polynomials.  We denote by $Q$
the categorical quotient $G \, \bb \, V \cong W \, \bb \, c$.
As a variety, $Q$ is isomorphic to a vector space of the same
dimension as $c$.

(v)   Let $f : V \to Q$ be the quotient map.  Denote by $c^{reg}
\subset c$ the set of regular points of $f|_c$.  It is the
complement of the union of the reflection hyperplanes in $c$.
Then the set $Q^{reg} \subset Q$ of regular values of $f$ is
equal to $f (c^{reg})$.
\end{thm}

We write $V^\circ = f^{-1} (Q^{reg})$.

\begin{rmk}\label{dkconj}
{\em  It is easy to see that $V^s \cap V^\circ \subset V^{rs}$.
Dadok and Kac conjecture that, in fact, $V^s \cap V^\circ = V^{rs}$
[DK, p. 521, Conjecture 2].}
\end{rmk}

We will call the fundamental group $\pi_1 (Q^{reg})$ the
{\em braid group} of $W$, and denote it by $B^{}_W$.

\begin{prop}\label{polardecom} {\em [DK]}
Let $G \, | \, V$ be a polar representation, as in Theorem
\ref{polar}, and $c \subset V$ a Cartan subspace.  Let
$G_c$ be the identity component of the stabilizer
$Z^{}_G (c)$.

(i)   There exists a compact form $K \subset G$, and a
$K$-invariant inner product $\langle \, , \, \rangle$ on $V$,
such that each $v \in c$ is of minimal length in its $G$-orbit. 

(ii)  The group $G_c$ is reductive, and $K_c = G_c \cap K$ is a 
compact form of $G_c$.

(iii) We have a $G_c$-invariant orthogonal decomposition
$$V = c \oplus U_c \oplus \fg \cdot c,$$
with $\dim \, G_c \bb U_c = 0$.

(iv)  The representation $G \, | \, V$ is stable if and only if
$U_c = 0$. 
\end{prop}

We let $d = \dim \, V$, $d_0 = d - \dim \, U_c$, and
$r = \dim \, c = \dim \, Q$.  Note that $V^{rs}$ is connected and
smooth of dimension $d_0$.

The rank of a polar representation $G \, | \, V$ is defined by
$$\rank \, G \, | \, V = \dim \, c / V^G,$$
where $V^G \subset V$ is the fixed points of $G$.  Note that
$V^G$ is contained in every Cartan subspace of $V$. 
Dadok and Kac [DK, Theorem 2.12] show that any polar
representation has a decomposition into polar representations
of rank one, analogous to the root space decomposition of a
semisimple Lie algebra.  We will return to this result in
Section 5.

\subsubsection{Geometry of the General Fiber of a Polar Quotient
Map}

We now deduce from the results of [DK] some corollaries about
the conormal geometry of the general fiber of a polar quotient
map.  First, we need to discuss the dual of a polar representation.
We continue with the situation of Theorem \ref{polar} and
Proposition \ref{polardecom}.

\begin{prop}\label{polardual}
(i)  The dual representation $G \, | \, V^*$ is also polar, with
the same Weyl group $W$.  Given a Cartan subspace $c \subset V$,
as in Proposition \ref{polardecom}, a Cartan subspace for
$G \, | \, V^*$ is given by
$$c^* = (\fg \cdot c \oplus U_c)^\perp \subset V^*.$$
This establishes a one-to-one correspondence between the
Cartan subspaces of $V$ and $V^*$.  Let
$\langle \, , \, \rangle^{}_{V^*}$ be the Hermitian metric on
$V^*$ induced by $\langle \, , \, \rangle$.  Then  every $l \in
c^*$ is of minimal length in its $G$-orbit.

(ii) Let $l \in c^*$ be regular semisimple.  Then the set
$\{ \, v \in V \; | \; l \, | \, _{\fg \cdot v} = 0 \, \}$
equals $c \oplus U_c$.
\end{prop}
\begin{pf}
The dual representation $G \, | \, V^*$ differs from $G \, | \, V$
by an outer automorphism of $G$.  Thus, it is also polar with a
Cartan subspace of the same dimension, and the same Weyl group.

Take any $l \in c^*$.  It follows from Proposition \ref{polardecom}
that $\langle \, \fg \cdot l, \, l \, \rangle^{}_{V^*} = 0$.  By a
theorem of Kempf and Ness ([KN], [DK, Theorem 1.1]), this implies
that $l$ is semisimple and of minimal length in its orbit.  The
rest of the proof is an exercise in linear algebra using
Proposition \ref{polardecom}.
\end{pf}

\begin{cor}\label{cortodual}
(i)   The set $(V^*)^{rs}$ of regular semisimple points in $V^*$ is
a connected algebraic manifold of dimension $d_0$.  We have
$(V^*)^s \subset \overline{(V^*)^{rs}}$.

(ii)  The subspace $U_c \subset V$ of Proposition \ref{polardecom}
is canonically independent of the choice of the compact form $K$,
and the inner product $\langle \, , \, \rangle$.
\end{cor}
\Qed

Pick a regular value $\lambda \in Q^{reg}$, and let
$F = f^{-1} (\lambda) \subset V$.  Let $\cO \subset F$ be the
unique closed $G$-orbit in $F$, and let $v \in c \cap \cO$.

\begin{prop}\label{polarfiber}
(i)   There is a unique $G$-equivariant algebraic map
$\phi : F \to \cO$, such that $\phi (v + u) = v$, for any
$u \in U_c$.

(ii)  The map $\phi : F \to \cO$ is a complex algebraic vector
bundle of rank $d - d_0$.

(iii) Fix any $x \in \cO$.  Then the tangent spaces $T_y F$,
for all $y \in \phi^{-1} (x)$, are parallel.
\end{prop}
\begin{pf}
Let $\psi : G \times U_c \to V$ be the map $\psi : (g, u)
\mapsto g \, (u + v)$.  By part (iii) of Proposition
\ref{polardecom}, we have $\mbox{\rm Im} \, \psi \subset F$.  It
is not hard to see that for small $u \in U_c$, the differential
$d_{(1, u)} \psi$ is surjective onto the tangent space
$T_{u+v} F$, which is parallel to $U_c \oplus \fg \cdot c$.
It follows that $\mbox{\rm Im} \, \psi$ contains a neighborhood
of $v$ in $F$.  Since $\cO$ is the unique closed $G$-orbit in
$F$, we have $\mbox{\rm Im} \, \psi = F$.  It also follows that
$T_{u+v} F$ is parallel to $U_c \oplus \fg \cdot c$, for all
$u \in U_c$.  It is now not hard to check that setting
$\phi (\psi (g, u)) = g \, v$ gives the required definition.
\end{pf}

We can now describe the conormal geometry of $F$.  Consider the
conormal bundle $\Lambda_F = T_F^* V \subset T^* V = V \times
V^*$, and let $p_2 : \Lambda_F \to V^*$ be the projection map.

\begin{cor}\label{polarcon}
(i)   The image $p_2 (\Lambda_F) = (V^*)^s$.

(ii)  For any $l \in (V^*)^{rs}$, the critical points of the
restriction $l |_F$ form a non-degenerate (Morse-Bott) critical
manifold, which is a union of $|W|$ fibers of $\phi$.   In
particular, $l |_F$ is a (complex) Morse function if and only
if $G \, | \, V$ is stable.
\end{cor}
\begin{pf}
Part (i) follows from the $G$-equivariance of the map $p_2$,
and Propositions \ref{polardual}, \ref{polarfiber}.  For part
(ii) we may assume that $l \in c^*$.  Then the critical points
of $l \, |_F$ is the set $\phi^{-1} (W \cdot v)$.  The
$G$-equivariance of $p_2$ and Proposition \ref{polardual} imply
that the restriction of $p_2$ to $p_2^{-1} (V^*)^{rs}$ is a
smooth submersion.  The Morse-Bott property of $l \, |_F$ is
simply a restatement of that fact.
\end{pf}

\begin{prop}\label{polarangle}
The variety $F \subset V$ is transverse to infinity in the sense
of Section 2.2.2, i.e., using the inner product of Proposition
\ref{polardecom}, we can find a $k > 0$, such that:
$$\angle (x, T_x X) < \frac{k}{|\!| x |\!|},$$
for any $x \in F$.
\end{prop}
\begin{pf}
We will show that we can take $k = |\!| \, v \, |\!|$.  Let
$x \in F$ and $l \in V^*$ be such that
$(x, l) \in \Lambda_F \subset T^* V \cong V \times V^*$.  Then
it is enough to show that $| l (x) | \leq 
|\!| \, v \, |\!| \cdot |\!| \, l \, |\!|$.  By Proposition
\ref{polarfiber}, we may write $x = g \, (v + u)$ and
$l = g \, l'$, where $g \in G$, $u \in U_c$, and $l' \in c^*$.
Then we have $l (x) = l' (v + u) = l' (v)$, and
$|l (x)| = |l' (v)| \leq |\!| \, v \, |\!| \cdot
|\!| \, l' \, |\!| \leq |\!| \, v \, |\!| \cdot
|\!| \, l \, |\!|$, where the last inequality follows from the
fact that $l'$ is of minimal length in its orbit (Proposition
\ref{polardual}).
\end{pf}

\section{The Main Theorem}

\subsection{Statement of the Theorem}

Let $G \, | \, V$ be a polar representation, as in Theorem
\ref{polar}.  Assume that either $\rank \, G \, | \, V = 1$,
or $G \, | \, V$ is {\em visible}, i.e., there are finitely
many nilpotent orbits in $V$.  This assumption holds for the
representations discussed in [KR] and for all of the infinite
series listed in [DK] and [K].  As before, we write $f : V \to Q$
for the quotient map, and $E = f^{-1} (0)$ for the zero fiber.
The purpose of the additional restriction of $G \, | \, V$
is to ensure that there exists a $G$-invariant conical
stratification $\cS$ of $E$, such that the ${\rm A}_f$
condition holds for the pair $(V^\circ, S)$, for every
$S \in \cS$.  Indeed, if $G \, | \, V$ is visible, we can
take $\cS$ to be the orbit stratification.  If, instead,
$\rank \, G \, | \, V = 1$, we can readily reduce to the case
$\dim \, Q = 1$, then use the general result about the
existence of ${\rm A}_f$ stratifications for functions
[Hi, p. 248, Corollary 1].

We are now in the situation of Proposition \ref{nclb}, and
therefore may consider the nearby cycles sheaf $P = P_f \in
\cP^{}_G \, (E)$.  In order to fix the up-to-isomorphism
ambiguity in the definition of $P$, we fix a regular value
$\lambda \in Q^{reg}$, let $F = f^{-1} (\lambda)$, and
identify $P$ with the sheaf $P^{}_F$ given by the specialization
of $F$ to $\as (F) = E$, as in Section 2.2.2.  This
corresponds to specializing along the path $\gamma : z \mapsto
f(z F)$.  We also have the monodromy action $\mu : B^{}_W
= \pi_1 (Q^{reg}, \lambda) \to {\rm Aut} (P)$.  Our main result
is Theorem \ref{main} below.  It gives a description of the
pair $(P, \mu)$ analogous to Theorem \ref{ftinspth} in Springer
theory.  Recall from Corollary \ref{cortodual} that the set
$(V^*)^{rs}$ of regular semisimple points in $V^*$ is a
connected algebraic manifold.

\begin{thm}\label{main}
Let $G \, | \, V$ be, as above, a polar representation which
is of rank one or visible.  Let $c \subset V$ be a Cartan
subspace, and $W$ be the Weyl group.

(i)    We have
$$\pF \, P \cong \IC ((V^*)^{rs}, \cL),$$
where the right hand side is an intersection homology sheaf
with coefficients in a local system $\cL$ on $(V^*)^{rs}$
of rank $|W|$.

(ii)   Write $A = {\rm End} (P)$, the endomorphisms of $P$
as a perverse sheaf (forgetting the $G$-equivariant structure).
We have $\dim \, A = |W|$.  The monodromy action $\mu$ gives
a surjection $\C [B^{}_W] \to A$.  

(iii)  Let $\sigma \in  W$ be a primitive reflection of order
$n_\sigma$.  It gives rise to an element $\hat\sigma \in B^{}_W$,
represented by a loop going counterclockwise around the image
$f(c_\sigma)$ of the hyperplane $c_\sigma \subset c$ fixed by
$\sigma$.  The minimal polynomial $R_\sigma$ of
$\mu (\hat\sigma) \in A$ has integer coefficients and is of
degree $n_\sigma$ (we normalize $R_\sigma$ to have leading
coefficient $1$).

(iv)   Fix a basepoint $l \in (V^*)^{rs}$.  By part (i), we have
an action
$$A = {\rm End} (P) \cong {\rm End} (\cL) \to {\rm End} (\cL_l)$$
of $A$ on the fiber $\cL_l$.  There is an identification
$\chi : \cL_l \cong A$, such that this action is given by left
multiplication.

(v)    There is a semigroup homomorphism $\rho : \pi_1
((V^*)^{rs}, l) \to A^{\circ}$, the opposite of the algebra $A$,
such that the holonomy of $\cL$ is given as $\rho$ followed by the
right multiplication action of $A^{\circ}$ on $\cL_l
\stackrel{\chi} {\cong} A$.  The homomorphism $\rho$ gives a
surjection $\C [\pi_1 ((V^*)^{rs}, l)] \to A^{\circ}$.
\end{thm}

\begin{rmk}\label{visibleispolar}
{\em Dadok and Kac conjecture [DK, p. 521, Conjecture 4] that
every visible representation is polar.}
\end{rmk}

\begin{rmk}\label{supp}
{\em  The support $\, \mbox{supp} \, (\pF \, P) =
\overline{(V^*)^{rs}} \,$ is equal to all of $V^*$ if and only if
the action $G \, | \, V$ is stable (cf. Corollary \ref{cortodual}).}
\end{rmk}

\begin{rmk}\label{stha}
{\em  It is natural to expect from claims (ii) and (iii) of
Theorem \ref{main} that the algebra $A$ is equal to the quotient
$\C [B^{}_W] / \langle \, R_\sigma (\hat\sigma) \, \rangle$, where
$\sigma$ runs over all the primitive reflections in $W$.  This is
always true when $W$ is a Coxeter group.  In fact, it is the case
in all of the examples known to the author.}
\end{rmk}

\begin{rmk}\label{pointer}
{\em  We will give some additional information about the
polynomials $R_\sigma$ in Section 5.}
\end{rmk}

\begin{rmk}\label{genft}
{\em  One may ask whether it is always true, say, for a
homogeneous polynomial $f : \C^d \to \C$ that the Fourier
transform $\pF \, P_f$ of the nearby cycles of $f$ is an
intersection homology sheaf.  The answer is ``no.''  A
counterexample is given by $f (x,y,z) = x^2 y + y^2 z$.}
\end{rmk}

\subsection{Examples}

\begin{ex}\label{qdrcs}
{\bf Quadrics.}
{\em Let $G = SO_n$ act on $V = \C^n$ by the standard
representation.  Any non-isotropic line $c \subset V$ can
serve as a Cartan subspace.  The Weyl group $W = \Z / 2 \Z$,
and the invariant map $f : V \to \C$ is just the standard
quadratic invariant of $SO_n$.  The algebra $A$ of Theorem
\ref{main} is given by $A = \C [z] / (z-1)^2$, if $n$ is even,
and $A = \C [z] / (z^2-1)$, if $n$ is odd.}
\end{ex}

\begin{ex}\label{nc}
{\bf Normal Crossings.}
{\em Consider the action of the torus $G = (\C^*)^{n-1}$ on
$V = \C^n$, given by
$$(t_1, \dots , t_{n-1}) : (x_1, \dots , x_n) \mapsto (t_1 x_1,
\, t_1^{-1} t_2 x_2, \, t_2^{-1} t_3 x_3, \, \dots , \,
t_{n-1}^{-1} x_n).$$
A Cartan subspace for this action is given by
$c = \{ \, x_1 = \dots = x_n \, \} \subset V$, the Weyl group
$W = \Z / n \Z$, and the invariant map $f : V \to \C$ is just the
product $f : (x_1, \dots , x_n) \mapsto x_1 \cdot \dots \cdot x_n$.
The algebra $A$ of Theorem \ref{main} is given by $A = \C [z] /
(z-1)^n$, functions on the $n$-th order neighborhood of a point.}
\end{ex}

\begin{ex}\label{det}
{\bf The Determinant} {\em (see [BG]).
Let $G = SL_n$ act on $V = Mat \, (n \times n; \C)$
by left multiplication.  A Cartan subspace $c \subset V$ for this
action is given by the scalar matrices, the Weyl group
$W = \Z / n \Z$, and the invariant map $f : V \to \C$ is the
determinant.  As in the previous example, the algebra
$A = \C [z] / (z-1)^n$.}
\end{ex}

\begin{ex}\label{symmat}
{\bf Symmetric Matrices} {\em (see [BG]).
Let $G = SL_n$ act on the space $V$ of symmetric $n\times n$
matrices by $g : x \mapsto g x g^t$.  A Cartan subspace
$c \subset V$ is again given by the scalar matrices, the Weyl group
$W = \Z / n \Z$, and the invariant map $f : V \to \C$ is again
the determinant.  However, the algebra $A$ is now given by
$A = \C [z] / (z-1)^{\lceil n/2 \rceil} (z+1)^{\lfloor n/2
\rfloor}$.}
\end{ex}

\begin{ex}\label{ran}
{\bf A ``Real Analog'' of Springer Theory for $SL_n$}
{\em (see [Gr2] and Section 6 below).
Let $G = SO_n$ act on the space $V$ of symmetric $n\times n$
matrices by conjugation.  A Cartan subspace $c \subset V$ is
given by the diagonal matrices.  The Weyl group $W = \Sigma_n$,
the symmetric group on $n$ letters.  The invariant map
$f : V \to \C^{n-1}$ is given by sending a matrix $x \in V$ to
the characteristic polynomial $\mbox{char} \, (x)$.  The algebra
$A = \cH^{}_{-1} (\Sigma_n)$, the Hecke algebra of $\Sigma_n$
specialized at $q = -1$.}
\end{ex}

\begin{ex}\label{fullstgp}
{\bf Maps from an Orthogonal to a Symplectic Vector Space.}
{\em Let $(U_1, \nu)$ be an orthogonal vector space of dimension
$2n+1$ ($\nu$ is a non-degenerate quadratic form), and
$(U_2, \omega)$ a symplectic vector space of dimension $2n$.  Set
$V = \mbox{Hom}^{}_\C \, (U_1, U_2)$.  The group $G = Sp_{2n}
\times SO_{2n+1}$ acts on $V$ by left-right multiplication.  This
action is polar, with a Cartan subspace of dimension $n$.  The Weyl
group $W \subset GL_n$ is the semi-direct product of $(\Z / 4 \Z)^n
\subset GL_n$, acting by diagonal matrices of fourth roots of unity,
and the symmetric group $\Sigma_n \subset GL_n$, acting by
permutation matrices.  Let $\sigma \in \Sigma_n \subset W$ be any
simple reflection, and let $\tau = \mbox{diag} \, (i, 1, \dots, 1)
\in (\Z / 4 \Z)^n \subset W$.  The algebra $A$ is given by
$A = \C [B^{}_W] / \langle (\sigma^2 - 1), \,
(\tau^2 - 1)^2 \rangle$.}
\end{ex}

In all of the above examples, the representation $G \, | \, V$ is
stable.  Here is an example where it is not.

\begin{ex}\label{nonstab}
{\bf An Evaluation Map.}
{\em  Let $U \cong \C^{2n}$.  Set $V = \Lambda^2 U^* \oplus U
\oplus U$.  The linear group $G = GL(U)$ acts on $V$.  The only
invariant is the evaluation map $f : (\omega, u_1, u_2) \mapsto
\omega (u_1, u_2)$.  A non-zero vector $(\omega, u_1, u_2) \in V$
is semisimple, if and only if $\omega$ is of rank $2$ and
$\omega (u_1, u_2) \neq 0$.  The Weyl group $W = \Z / 3 \Z$, and
the algebra $A = \C [z] / (z - 1)^3$.}
\end{ex}

\subsection{A Remark on the Degree of Generality}

The reason we assume in Theorem \ref{main} that $G \, | \, V$ is
of rank one or visible, is to insure that the ${\rm A}_f$ condition
holds and the sheaf $P$ is well defined.  If we dropped this
assumption, we could still pick a $\lambda \in Q^{reg}$, and
consider a sheaf $P_\lambda$ coming from the specialization of
$F = f^{-1} (\lambda)$ to the asymptotic cone.  Using the techniques
of this paper and some algebro-geometric generalities, we could then
show that the sheaves $P_\lambda$ for {\em generic} $\lambda$ are
isomorphic, and that they form a local system over some open
set $Q' \subset Q^{reg}$.  Furthermore, all of the assertions
of Theorem \ref{main} could be given meaning and proved in this
context.  This is somewhat unsatisfactory, since we would
have to use the full power of our methods just to show that there
is a well posed question.  Conjecture \ref{polaraf} below implies
that the subset $Q'$ is a phantom, and that nearby cycles are well
defined for any polar representation.

\begin{conj}\label{polaraf}
Let $G \, | \, V$ be a polar representation, $f : V \to Q$ the
quotient map, and $E = f^{-1} (0)$.  Then there exists an
${\rm A}_f$ stratification of $E$ (see Section 2.2.1).
\end{conj}

If this conjecture is true, all our results extend automatically
to an arbitrary polar representations.

\section{Proof of the Main Theorem}

\subsection{A Preliminary Lemma}

We begin with a preliminary lemma.  Recall the $G$-invariant
conical stratification $\cS$ of $E$ discussed at the beginning
of Section 3.  By passing if necessary to a refinement, we
may assume that $\cS$ satisfies conditions (i) - (iii) of
Section 2.2.2.

\begin{lemma}\label{prelim}
Let $l \in (V^*)^{rs}$, and let $a (\cS, l)$ be the dimension
of the critical locus of $l \, |_E$ with respect to
$\cS$, as in Lemma \ref{dimbound}.  Then $a (\cS, l) \leq
d - d_0$, where $d = \dim V$, and $d_0 = \dim \, V^{rs}$.
\end{lemma}
\begin{pf}
Let $\Lambda_E \subset T^* V = V \times V^*$ be the conormal
variety to the stratification $\cS$, and $p_2 : V \times V^*
\to V^*$ be the projection.  Using Proposition \ref{polardual}
and the $G$-invariance of $\cS$, it is not hard to show that
$\dim \, p_2^{-1} (l) \cap \Lambda_E$ is independent of $l$
for $l \in (V^*)^{rs}$.  Therefore,
$$a (\cS, l) + \dim \, (V^*)^{rs} \leq \dim \, \Lambda_E = d.$$
Together with Corollary \ref{cortodual}, this proves the lemma.
\end{pf}

\subsection{The Stable Case}

In this section we prove Theorem \ref{main} in the case when
the representation $G \, | \, V$ is stable.  Then, in Section
4.3, we will indicate how to modify the argument in the
nonstable case.

Assuming $G \, | \, V$ is stable, choose a basepoint
$l \in (V^*)^{rs}$.  Note that $(V^*)^{rs}$ is open in $V^*$.
We may assume that $l \in c^*$, the Cartan subspace of
Proposition \ref{polardual}.  By Proposition \ref{polarangle},
all the constructions of Section 2.2.2 apply to $X = F$.  As
in that section, we consider the compactification $\hat F$ of
$F$ relative to $l$, and denote by $Z$ the set of critical
points of the restriction $\hat l \, |_{\hat F}$, as in Lemma
\ref{dimbound}.  By Lemmas \ref{dimbound}, \ref{prelim}, and
Corollary \ref{polarcon}, we have $Z \subset c \subset V$.
Furthermore, $Z$ is a single $W$-orbit in $c$, and each
critical point in $Z$ is Morse. Fix a point $e^{}_0 \in Z$,
and write $Z = \{ e^{}_w \}^{}_{w \in W}$, where $e^{}_w =
w \, e^{}_0$.

Recall the description of the stalk $H^{*}_l \, (\pF \, P)$ in
Lemma \ref{auxmtw}.  Since $\dim \, Z = 0$, we have
$H^{i}_l \, (\pF \, P) = 0$, unless $i = -d$.  Furthermore,
by duality:
$$H^{-d}_l \, (\pF \, P) \cong H_{d-r} \, ( F, \, \{ \, \xi (y)
\geq \xi^{}_0 \, \}; \; \C),$$
where $\xi^{}_0 > | \, l (e) \, |$, for all $e \in Z$.

We will need a standard Picard-Lefschetz construction
of classes in the relative homology group above (see,
for example, [BG, Section 7.2]).  Let $e \in Z$ be a critical
point.  Fix a smooth path $\gamma : [0, 1] \to \C$ such that:

(i)   $\gamma (0) = l (e)$, and $\gamma (1) = \xi^{}_0$;

(ii)  $\gamma (t) \notin l (Z)$, for $t > 0$;

(iii) $\gamma (t_1) \neq \gamma (t_2)$, for $t_1 \neq t_2$;

(iv)  $\gamma' (t) \neq 0$, for $t \in [0, 1]$.

\noindent
Let $\cH_e : T_e F \rightarrow \C$ be the Hessian of
$l \, |_F$ at $e$, and let $T_e \, [\gamma] \subset T_e F$ be the
positive eigenspace of the (non-degenerate) real quadratic form
$$\mbox{Re} \, (\, {\cal H}_e \, / \, \gamma'(0) \, ) :
T_e F \to \R.$$
Note that $\dim_\R \, T_e \, [\gamma] = d-r$.  Fix an orientation
${\cal O}$ of $T_e \, [\gamma]$.  The triple $(e, \gamma, \cO)$
defines a homology class
$$u = u (e, \gamma, \cO) \in H_{d-r} \, ( F, \,
\{ \, \xi (y) \geq \xi^{}_0 \, \}; \; \C).$$
Namely, the class $u$ is represented by an embedded $(d-r)$-disc
$$\kappa : (\D^{d-r}, \, \partial \, \D^{d-r}) \to
(F, \, \{ \, \xi (y) \geq \xi^{}_0 \, \}),$$
such that the image of $\kappa$ projects onto the image of
$\gamma$, is tangent to $T_e \, [\gamma]$ at $e$, and does not
pass through any point of $Z$ except $e$.  The sign of $u$ is given
by the orientation $\cO$.  It is a standard fact that $u \neq 0$.

We now fix a path $\gamma^{}_0 : [0, 1] \to \C$, satisfying
conditions (i) - (iv) above for $e = e^{}_0$, and an orientation
$\cO^{}_0$ of the space $T_{e^{}_0} \, [\gamma^{}_0]$.  Let
$u^{}_0 = u (e^{}_0, \gamma^{}_0, \cO^{}_0)$; we will use the
same symbol for the corresponding element of $H^{-d}_l \,
(\pF \, P)$.

\begin{lemma}\label{monsurj}
The image of $u^{}_0 \in H^{-d}_l \, (\pF \, P)$ under the
monodromy action of $B^{}_W$ generates the stalk
$H^{-d}_l \, (\pF \, P)$ as a vector space.  We have
$\dim \, H^{-d}_l \, (\pF \, P) = |W|$.
\end{lemma}
\begin{pf}
For each $w \in W$, pick a lift $b^{}_w \in B^{}_W$.  By a
standard Picard-Lefschetz theory argument, we have
$$\mu (b^{}_w) \, u^{}_0 = u \, (e^{}_w, \gamma^{}_w, \cO^{}_w),$$
where $\gamma^{}_w$ is some path satisfying conditions (i) - (iv)
for $e = e^{}_w$, and $\cO^{}_w$ is an orientation of
$T_{e^{}_w} \, [\gamma^{}_w]$.  Since $\{ e^{}_w \}^{}_{w \in W}$
are the only critical points of the function $\hat l \, |_{\hat F}$,
the elements $\{ \mu (b^{}_w) \, u^{}_0 \}^{}_{w \in W}$ form a
basis of $H^{-d}_l \, (\pF \, P)$.
\end{pf}

Consider now the restriction of $\pF \, P$ to the manifold
$(V^*)^{rs}$.  By Lemma \ref{monsurj} and the preceding arguments,
it is a perverse sheaf whose stalks only live in degree $-d$ and
all have rank $|W|$.  Therefore, it is a rank $|W|$ local system
with a degree shift.  Together with Theorem \ref{asym_cone}, this
proves part (i) of Theorem \ref{main}.  The local system $\cL$ of
that theorem is given by $\cL = \pF \, P \, |_{\, (V^*)^{rs}} \,
[-d]$.  Let $h :  \pi_1 ((V^*)^{rs}, l) \to {\rm End} (\cL_l)$ be
the holonomy of $\cL$.

\begin{lemma}\label{holsurj}
The image of $u^{}_0 \in \cL_l$ under the action $h$ generates
the stalk $\cL_l$ as a vector space.
\end{lemma}
\begin{pf}
The proof of this is analogous to the proof of Lemma
\ref{monsurj}.  Let $Q' = G \bb V^*$, and $f' : V^* \to Q'$
be the quotient map.  Let $(Q')^{reg} \subset Q'$ be the set
of regular values of $f'$.   Set $(V^*)^\circ = {f'}^{\, -1}\,
((Q')^{reg})$; it is a Zariski open subset of $(V^*)^{rs}$ (see
Remark \ref{dkconj}).  We may assume without loss of generality
that $l \in (V^*)^\circ$.  Let $\cL^\circ_{}$ be the restriction
of the local system $\cL$ to $(V^*)^\circ$, and let $h^\circ :
\pi_1 ((V^*)^\circ, l) \to {\rm End} (\cL_l)$ be the holonomy of
$\cL^\circ$.  It is enough to show that the image of $u^{}_0$
under $h^\circ$ generates the stalk $\cL_l$.

Since we assume that $G$ is connected, the push-forward
homomorphism $f'_* : \pi_1 ((V^*)^\circ, l) \to
\pi_1 ((Q')^{reg}, f' (l))$ is surjective.  By Proposition
\ref{polardual}, we have an isomorphism
$\eta : \pi_1 ((Q')^{reg}, f' (l)) \cong B^{}_W$ (it is not
canonical).  The composition $\eta \circ f'_*$, followed by the
natural map $B^{}_W \to W$, gives a surjection
$\theta : \pi_1 ((V^*)^\circ, l) \to W$.  For each $w \in W$,
choose a lift $\alpha^{}_w \in \theta^{-1} (w)$.  Then a
Picard-Lefschetz theory argument similar to the proof of Lemma
\ref{monsurj} shows that the elements
$\{ \, h^\circ (\alpha^{}_w) \, u^{}_0 \, \}^{}_{w \in W}$
form a basis of $\cL_l$.
\end{pf}

Claims (ii), (iv), and (v) of Theorem 3.1 follow easily from
Lemmas \ref{monsurj} and \ref{holsurj}.  We define the vector 
space map $\chi : A \to \cL_l$ of claim (iv) by $\chi : a \mapsto
a \, u^{}_0$.  By Lemma \ref{monsurj}, $\chi$ is surjective.  On
the other hand, the action of any $a \in A$ on $\cL_l$ must
commute with the holonomy representation $h$.  Therefore, Lemma
\ref{holsurj} implies that $\chi$ is injective.  This proves that
$\chi$ is a vector space isomorphism.  Claim (ii) now follows from
Lemma \ref{monsurj}; claim (iv) from the definition of $\chi$; and
claim (v) from Lemma \ref{holsurj}, and the fact that every
endomorphisms of $A$ commuting with the left action of $A$ on
itself is of the form $a \mapsto a \, a'$, for some $a' \in A$
(this is true for any associative algebra with unit).

It remains to prove claim (iii) of the theorem.  The idea is to
use a regular value $\lambda \in Q^{reg}$ which is very near the
image $f (c_\sigma)$, then apply a standard carousel argument
(see [L\^e2] and [BG, Section 7.3]).

Pick a point $v_1 \in c_\sigma \subset c$ which is not fixed
by any element of $W$, other than the powers of $\sigma$.  The
orbit $W \cdot v_1$ consists of $|W| / n_\sigma$ points.  Choose
a point $v \in c^{reg}$ which is very near $v_1$.  The orbit
$W \cdot v$ consists of $|W|$ points which are grouped into
$|W| / n_\sigma$ clusters.  Each cluster consists of $n_\sigma$
points surrounding a point in $W \cdot v_1$, and the action of
$\sigma \in W$ cyclically permutes the points in each cluster. 
We may assume that the fiber $F = f^{-1} (\lambda)$ passes
through $v$, so that $Z = W \cdot v$.  We may also assume that
the choice of $l \in c^* \cap (V^*)^{rs}$ is sufficiently generic,
so that the image $l (Z)$ appears in the $\C$-plane as
$|W| / n_\sigma$ disjoint clusters, each consisting of $n_\sigma$
points arranged in a circle (see Figure 1).

Choose a path $\gamma_0 : [0, 1] \to \C$, satisfying conditions
(i) - (iv) above for $e = v$, which does not ``come near'' any
of the clusters in $l (Z)$, other than the one containing
$l (v)$.  Fix an orientation $\cO_0$ of the space
$T_{v} \, [\gamma]$.  Let $u_0 = u (v, \gamma_0, \cO_0)$.
The element $\hat \sigma \in B^{}_W = \pi_1 (Q^{reg}, \lambda)$ is
represented by a loop going once counterclockwise around the image
$f (c_\sigma)$, which stays in a small neighborhood of the critical
value $f (v_1)$.  Then, by a standard carousel argument, the vectors
$\{ \, \mu(\hat \sigma)^k \, u_0 \, \}_{k = 0}^{n_\sigma - 1}$
are linearly independent, and we have
$$\mu(\hat \sigma)^{n_\sigma} \, u_0 = \pm u_0 +
\sum_{k = 1}^{n_\sigma - 1} g_k \cdot \mu(\hat \sigma)^k \, u_0,$$
where the $\{ g_k \}$ are some integers.  Together with Lemma
\ref{holsurj}, this implies claim (iii) of Theorem \ref{main}.

\begin{center}
\leavevmode{}
\epsfbox{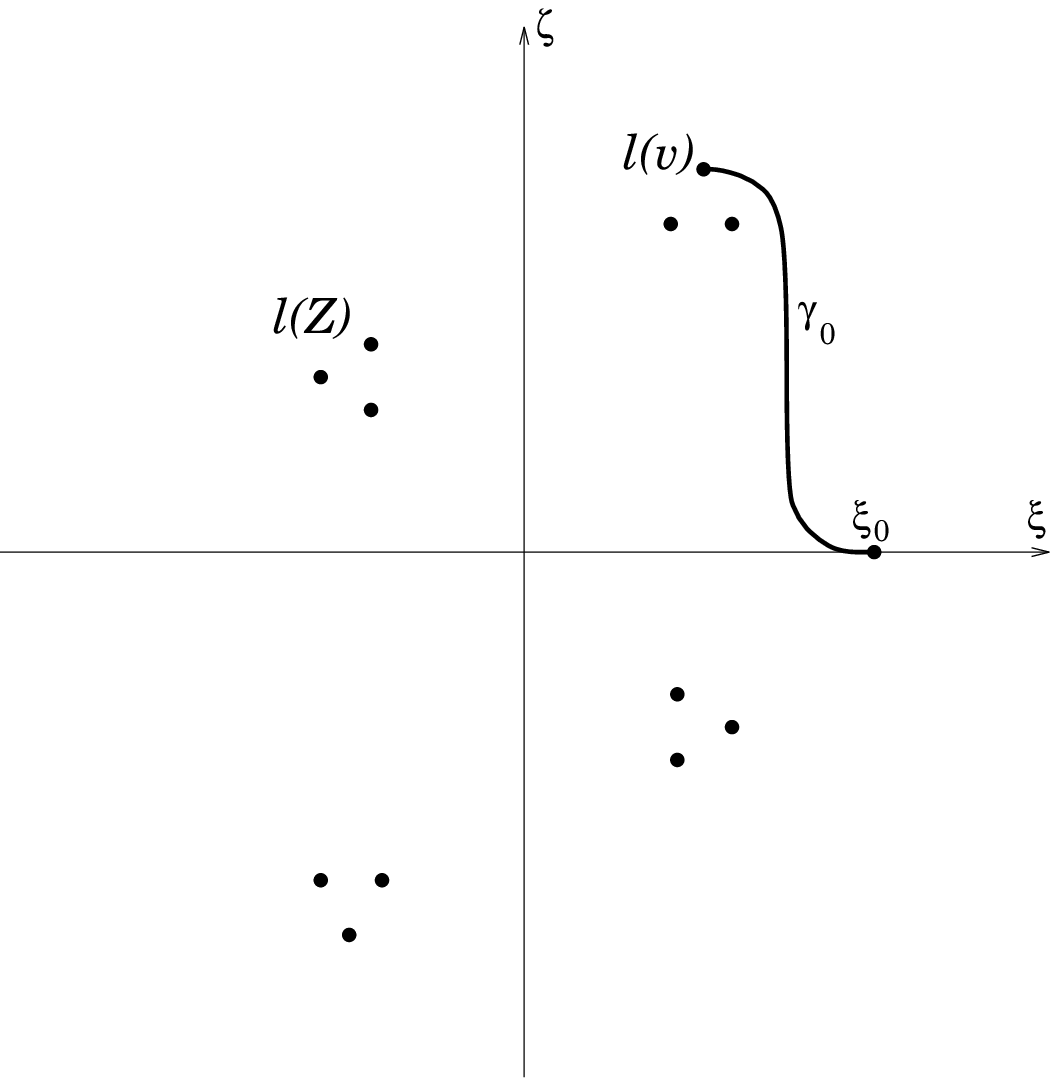}
Figure 1
\end{center}
\vspace{.1in}

\subsection{The Nonstable Case}

The proof of Theorem \ref{main} in the case when $G \, | \, V$ is
not stable is analogous to the argument of Section 4.2, with Morse
critical points of $l$ replaced by Morse-Bott critical manifolds
(see part (ii) of Corollary \ref{polarcon}).  We briefly indicate
the changes that have to be made in the nonstable case.

As before, we pick a basepoint $l \in (V^*)^{rs}$, lying in the
Cartan subspace $c^*$ of Proposition \ref{polardual}, and write $Z$
for the critical points of the restriction $\hat l \, |_{\hat F}$.
Fix a point $e \in (F \cap c) \subset Z$, and a path
$\gamma: [0,1] \to \C$, satisfying conditions (i) - (iv) of Section
4.2.  The Hessian $\cH_e : T_e F \rightarrow \C$ of $l \, |_F$ at
$e$ will now be degenerate.  However, we may still consider the
positive eigenspace $T_e \, [\gamma] \subset T_e F$ of the
real quadratic form
$\mbox{Re} \, (\, {\cal H}_e \, / \, \gamma'(0) \, )$.
By Corollary \ref{polarcon}, we have $\dim_\R \, T_e \, [\gamma] =
d_0 - r$.  Fix an orientation ${\cal O}$ of $T_e \, [\gamma]$.  The
triple $(e, \gamma, \cO)$ defines a Picard-Lefschetz class
$$u = u (e, \gamma, \cO) \in H_{d_0 - r} \, ( F, \,
\{ \, \xi (y) \geq \xi^{}_0 \, \}; \; \C),$$
where $\xi_0$ is large, exactly as before.  We may regard $u$ as an
element of $H^{-d_0}_l \, (\pF \, P)$.

The main distinction from the stable case is that there is no
a priori geometric reason for $u$ to be non-zero.  This is because
the Morse-Bott manifold containing $e$ is non-compact.  However,
we may use Lemmas \ref{auxmtw}, \ref{dimbound}, \ref{prelim}, and
an argument as in Lemma \ref{monsurj}, to show that $u$ generates
the stalk $H^{-d_0}_l \, (\pF \, P)$ under the monodromy action of
$B^{}_W$.  Thus, if we had $u = 0$, it would follow that
$H^{-d_0}_l \, (\pF \, P) = 0$, and by Theorem \ref{asym_cone},
that $P = 0$.  This contradiction shows that $u \neq 0$.  The rest of
the proof is exactly as in the stable case.

\section{Rank One Representations and the Polynomials
$R_\sigma$}

As we mentioned in Section 2, there is an analog of root space
decomposition for polar representation.  Namely, any polar
representation has a decomposition into polar representations
of rank one.  In this section, we use these rank one
representations to give an interpretation (see Theorem
\ref{rankone}) of the polynomials $R_\sigma$ appearing in part
(iii) of Theorem \ref{main}.  This interpretation will be used
in Section 6.2 to analyze the polar representations arising from
symmetric spaces.

Let $G \, | \, V$ be a polar representation, $\sigma \in  W$ be
a primitive reflection of order $n_\sigma$, and $c_\sigma
\subset c$ be the hyperplane fixed by $\sigma$ (cf. Theorem
\ref{main}).  Denote by $c_\sigma^\circ$ the set of all
$v \in c_\sigma$ such that the stabilizer $Z^{}_W (v)$ is
generated by $\sigma$; it is an open subset of the hyperplane
$c_\sigma$.  Let $\fg_\sigma \subset \fg$ be the stabilizer of
$c_\sigma$, and $G_\sigma \subset G$ the adjoint form of
$\fg_\sigma$.  In terms of Proposition \ref{polardecom}, let
$$V_\sigma = c \oplus U_c \oplus \fg_\sigma \cdot c \subset V.$$ 
If $\fg_\sigma \cdot c = 0$, we say that $\sigma$ is {\em of global
type}.  Otherwise, we say that $\sigma$ is {\em of local type}.

\begin{prop}\label{polarslice} {\em [DK]}
(i)   The representation $G \, | \, V$ restricts to a
representation $G_\sigma \, | \, V_\sigma$, which is polar
with Cartan subspace $c$.  The rank of $G_\sigma \, | \, V_\sigma$
is equal to zero if $\sigma$ is of global type, and to one if
$\sigma$ is of local type.

(ii)  Choose a point $v_1 \in c_\sigma^\circ$.  Let $\fg_{v^{}_1}
\subset \fg$ be the stabilizer of $v_1$.  Then we have
$\fg_{v^{}_1} = \fg_\sigma$, $\fg \cdot v_1 = \fg \cdot c_\sigma$,
and  $V = V_\sigma \oplus \fg \cdot v_1$.
\end{prop}

Proposition \ref{polarslice} implies that $\sigma$ is of
local type if and only if $c_\sigma \cap V^{rs} = \emptyset$.
Dadok and Kac conjecture (cf. Remark \ref{dkconj}) that $\sigma$ is
always of local type.  We do not know of any counterexamples, but
allow for the possibility that they exist.  The ``root space
decomposition'' theorem of [DK] involves the representations
$G_\sigma \, | \, V_\sigma$, for all primitive reflections $\sigma$
of local type.

The Weyl group $W_\sigma$ of $G_\sigma \, | \, V_\sigma$ is
naturally a subgroup of $W$.  Moreover, it must be a cyclic group,
generated by $\sigma^{m_\sigma}$ for some $m_\sigma \in \Z_+$, a
divisor of $n_\sigma$.  Note that when $\sigma$ is of global type,
we have $m_\sigma = n_\sigma$ and $W_\sigma = \{ 1 \}$.  Let
$f_\sigma : V_\sigma \to Q_\sigma = G_\sigma \bb V_\sigma$ be the
quotient map, $E_\sigma = f_\sigma^{-1} (0)$ be the zero fiber,
and $P_\sigma \in \cP_{G_\sigma} (E_\sigma)$ be the nearby cycles
of $f_\sigma$.  Let $B_\sigma$ be the braid group of $W_\sigma$.
If $\sigma$ is of local type, then $B_\sigma \cong \Z$, and we have
a monodromy transformation  $\mu_\sigma : P_\sigma \to P_\sigma$,
given by the action of the counterclockwise generator of $B_\sigma$.
If $\sigma$ is of global type, then $B_\sigma = \{ 1 \}$, and we set
$\mu_\sigma = \mbox{Id} : P_\sigma \to P_\sigma$.  Let $\tilde R$ be
the minimal polynomial of $\mu_\sigma \in {\rm End} (P_\sigma)$.  It
is a monic polynomial with integer coefficients of degree
$n_\sigma / m_\sigma$.  Note that if $\sigma$ is of global type,
then $\tilde R_\sigma (z) = z - 1$.  By Theorem \ref{main}, we have
${\rm End} (P_\sigma) \cong \C [z] / \tilde R_\sigma (z)$.

\begin{thm}\label{rankone}
Let $G \, | \, V$ be as in Theorem \ref{main}, and $G_\sigma
\, | \, V_\sigma$ be as in Proposition \ref{polarslice}.  Then the
polynomial $R_\sigma$ of part (iii) of Theorem \ref{main} is given
by $R_\sigma (z) = \tilde R_\sigma (z^{m_\sigma})$.
\end{thm}
\begin{pf}
Choose a point $v_1 \in c_\sigma$, as in part (ii) of Proposition
\ref{polarslice}.  Let $\lambda_1 = f (v_1) \in Q$.  Choose a
regular value $\lambda \in Q^{reg}$ very near $\lambda_1$.  The
idea of this proof is to break up the specialization of $\lambda$
to $0$ into two steps: first specialize $\lambda$ to $\lambda_1$,
then specialize $\lambda_1$ to $0$.  In order to compute the
polynomial $R_\sigma$, we will only need to understand the first
step.

Let $\D = \{ \, z \in \C \; | \; |z| < 2 \, \}$ and $\gamma^{}_1 :
\D \to Q$ be an embedded holomorphic arc such that:

(i)   $\gamma^{}_1 (0) = \lambda_1$;

(ii)  $\gamma^{}_1 (1) = \lambda$;

(iii) $\gamma^{}_1 (\D \setminus 0) \subset Q^{reg}$;

(iv)  $\gamma^{}_1$ is transverse to the image $f (c_\sigma)$.

\noindent
Form the fiber product $V_1 = V \times^{}_Q \D$, and let $f_1 : V_1
\to \D$ be the projection map.  Let $E_1 = f^{-1} (\lambda_1)$.  We
have $\dim \, E_1 = \dim \, E = d - r$.  Let $P_1 = \psi_{f_1} \,
\C^{}_{V_1} \, [d - r]$; we have $P_1 \in \cP^{}_G (E_1)$.  Write
$\mu_1 : P_1 \to P_1$ for the monodromy transformation of $P_1$.

{\bf Claim:}
There is a functor $\Psi : \cP_G (E_1) \to \cP_G (E)$ such
that $\Psi \, P_1 = P$ and $\Psi (\mu_1) = \mu (\hat\sigma)$.

To prove the claim, define an arc $\gamma_2 : \D \to Q$ by
$\gamma_2 : z \mapsto f (z \, v_1)$.  Form $V_2 = V \times^{}_Q \D$,
and let $f_2 : V_2 \to \D$ be the projection.  We have a map
$\pi : V_2 \setminus E \to E_1$, defined by $\pi (z \, e, \, z) = e$,
for $z \in \D \setminus 0$, $e \in E_1$.  Set $\Psi = \psi_{f_2}
\circ \pi^*$.  The claim follows form the definition of nearby
cycles.

Because of the claim, it will suffice to show that
$\tilde R_\sigma (\mu_1^{m_\sigma}) = 0$.  For this, consider
the inclusion $j_\sigma : V_\sigma \to V$ given by
$j_\sigma : v \mapsto v + v_1$.  By Proposition \ref{polarslice},
$j_\sigma$ exhibits $V_\sigma$ as a $G_\sigma$-invariant normal
slice to the orbit $G \cdot v_1$.  Consider the perverse
restriction to the normal slice functor $\pj_\sigma^*  =
j_\sigma^* \, [d_\sigma - d] : \cP_G (E_1) \to \cP_{G_\sigma}
(E_\sigma)$, where $d_\sigma = \dim \, V_\sigma$.  It is not
hard to see that $\pj_\sigma^*$ is injective on morphisms.  Set
$P_2 = \pj_\sigma^* \, P_1$, and $\mu_2 = \pj_\sigma^* (\mu_1) :
P_2 \to P_2$.  Then it will suffice to show that
\begin{equation}\label{EE}
\tilde R_\sigma (\mu_2^{m_\sigma}) = 0.
\end{equation}

Let $g_\sigma : Q_\sigma \to Q$ be the natural map.   Note
that $j_\sigma \circ f = f_\sigma \circ g_\sigma$.  Note also
that near the point $f_\sigma (v_1)$, the map $g_\sigma$ is an
$m_\sigma$-fold cover, ramified along the smooth hypersurface
$f_\sigma (c_\sigma)$.  It follows that
$P_2 \cong \oplus_{k = 1}^{m_\sigma} P_\sigma$,
and that $\mu_2$ is given by the matrix below:
$$
\pmatrix{0      & 0 & \dots  & 0      & \mu_\sigma \cr
         1      & 0 &        & 0      & 0          \cr
         0      & 1 &        & 0      & 0          \cr 
         \vdots &   & \ddots &        & \vdots     \cr
         0      & 0 & \dots  & 1      & 0          \cr}
$$
This implies (\ref{EE}).
\end{pf}

\begin{rmk}\label{partres}
{\em Our technique of first specializing $\lambda$ to $\lambda_1$,
then specializing $\lambda_1$ to $0$ has a counterpart in the
resolution approach to Springer theory.  Namely, it corresponds
to Borho and MacPherson's idea of first ``forgetting the
complete flag partially,'' then ``forgetting the partial flag
completely'' (see [BM2, p. 28]).}
\end{rmk}

\begin{cor}
All the roots of $R_\sigma$ are roots of unity.
\end{cor}
\begin{pf}
This follows form Theorem \ref{rankone} and the quasiunipotence
of the monodromy transformation $\mu_\sigma$ (see, for example,
[L\^e2]).
\end{pf}

\section{Symmetric Spaces}

Both Springer theory and Example \ref{ran} above are special
cases of representations arising from symmetric spaces studied
by Kostant and Rallis in [KR].  Let $\fg$ be a semisimple Lie
algebra over $\C$, and $\theta : \fg \to \fg$ an involutive
automorphism.  Let $\fg = \fg^+ \oplus \fg^-$ be the eigenspace
decomposition for $\theta$, so that $\theta |_{\fg^\pm} = \pm 1$.
Then $\fg^+$ is a Lie algebra, and the adjoint form $G^+$ of
$\fg^+$ acts on the symmetric space $\fg^-$ by conjugation.  By
the results of [KR], the representation $G^+ \, | \, \fg^-$ is
polar and visible.  Theorem \ref{symsp} below gives a recipe for
computing the algebra $A$ of Theorem \ref{main} for this
representation.

Let $\fg^{}_\R$ be a real form of $\fg$ with a Cartan
decomposition $\fg^{}_\R = \fg^+_\R \oplus \fg^-_\R$,
such that $\fg^\pm$ is the complexification of
$\fg^\pm_\R$ (see [OV] for a discussion of real forms
of complex semisimple Lie algebras).  Let
$c^{}_\R \subset \fg^-_\R$ be a maximal abelian
subalgebra.  Then the complexification
$c \subset \fg^-$ of $c^{}_\R$ is a Cartan subspace
for $G^+ \, | \, \fg^-$.  The Weyl group $W$ of this
representation is just the small Weyl group associated to
$\fg^{}_\R$.  It is a Coxeter group acting on $c^{}_\R$
with its Euclidean structure induced by the Killing
form on $\fg^{}_\R$.

Fix a Weyl chamber $C \subset c^{}_\R$.  Let
$\{ \, \sigma_i \, \}_{i=1}^r$ ($r = \dim^{}_\R \, c^{}_\R$)
be the reflections in the walls of $C$; they give a set
of generators for $W$.  Fix a basepoint $b \in C$.  Consider
the braid group $B^{}_W = \pi_1 (Q^{reg}, f(b))$ (we use
the notation of Theorem \ref{polar}).  Let $\hat\sigma_i \in
B^{}_W$ be the element represented by the $f$-image of a
path from $b$ to $\sigma_i (b)$ in $c^{reg}$, which is almost
straight but passes counterclockwise half-way around the
hyperplane $c_{\sigma_i} \subset c$ fixed by $\sigma_i$.  The 
braid group $B^{}_W$ is generated by the elements
$\{ \, \hat\sigma_i \, \}_{i=1}^r$.

To each reflection $\sigma_i$ we associate a number $s(i)$
as follows.  Let $\fg^+_c \subset \fg^+$ be the stabilizer
of $c$, and let $\fg^+_i \subset \fg^+$ be the stabilizer
of $c_{\sigma_i}$.  Then we set
$$s(i) = \dim \, \fg^+_i - \dim \, \fg^+_c.$$
(This is just half the sum of the dimensions of the root
spaces in $\fg^{}_\R$ corresponding to $\sigma_i$.)

The Killing form on $\fg$ restricts to a non-degenerate
$G^+$-invariant quadratic form on $\fg^-$.  We will use it
to identify $\fg^-$ with $(\fg^-)^*$, and to regard the
Fourier transform on $\fg^-$ as a functor
$\pF : \cP^{}_{\C^*} (\fg^-) \to \cP^{}_{\C^*} (\fg^-)$.
The set of regular semisimple vectors in $\fg^-$ is given by
$(\fg^-)^{rs} = f^{-1} (Q^{reg})$.

\begin{thm}\label{symsp}
Let $G^+ \, | \, \fg^-$ be the representation arising from an
involution $\theta$ as above.  Then, in terms of Theorem
\ref{main}, we have:

(i)   The endomorphism algebra $A = {\rm End} (P)$ is given as
\begin{equation}\label{FF}
A = \C [B^{}_W] / \langle \, (\hat\sigma_i - 1)
(\hat\sigma_i + (-1)^{s(i)}) \, \rangle_{i=1}^r.
\end{equation}

(ii)  Take $l = b$ as a basepoint for $(\fg^-)^{rs}$.
The homomorphism $\rho : \pi_1 ((\fg^-)^{rs}, \, b) \to
A^\circ$ of part (v) of Theorem \ref{main} is given by
$\rho = \alpha \circ f_*$, where
$f_* : \pi_1 ((\fg^-)^{rs}, \, b) \to \pi_1 (Q^{reg}, f(b)) =
B^{}_W$ is the push-forward by $f$, and $\alpha : B^{}_W \to
A^\circ$ is defined by $\alpha : \hat\sigma_i \mapsto
(-1)^{s(i) - 1} \cdot \hat\sigma_i$.
\end{thm}
\begin{pf}
Because $W$ is a Coxeter group, the braid group quotient in the
right hand side of (\ref{FF}) has the correct dimension, namely
the order of $W$.  Thus, for part (i) it is enough to show that
for $\sigma = \sigma_i$, we have $R_\sigma (z) = (z - 1)
(z + (-1)^{s(i)})$.  This is an easy application of Theorem
\ref{rankone}.  It follows from the root space decomposition for
$\fg^{}_\R$ that, in terms of Proposition \ref{polarslice}, we
have $G_\sigma = SO (s(i) + 1)$, $V_\sigma = \C^{s(i) + r}$, and
$G_\sigma$ acts by the standard representation on the first
$s(i) + 1$ coordinates.  The computation of $R_\sigma$ is now
given by Theorem \ref{rankone} and Example \ref{qdrcs}.

Part (ii) of the theorem is proved by a Picard-Lefschetz argument
as in Section 5 (cf. Lemma \ref{holsurj}).
\end{pf}

\begin{rmk}\label{hecvsgrp}
{\em  Note that $A$ is isomorphic to the group algebra
$\C [W]$, when $s(i)$ is even for all $i$, and to the
Hecke algebra $\cH^{}_{-1} (W)$ specialized at $q = -1$,
when $s(i)$ is odd for all $i$.  In general, it is a ``hybrid''
of the two.  Group algebra deformations of this kind were
considered by Lusztig in [Lu2].}
\end{rmk}

\vspace{.1in}

\noindent
MIT, Department of Mathematics, 77 Massachusetts Ave., Room 2-247,
Cambridge, MA 02139

\vspace{.1in}

\noindent
{\it grinberg@math.mit.edu}

\end{document}